\documentstyle[11pt]{article}
\newtheorem{theorem}{Theorem}
\newtheorem{proposition}{Proposition}
\newtheorem{lemma}{Lemma}
\newtheorem{corollary}{Corollary}

\def\IR{\Bbb R}
\def\tF{{\tilde F}}
\def\tG{{\tilde G}}
\def\tE{{\tilde E}}
\def\hF{{\hat F}}
\def\hG{{\hat G}}
\def\hE{{\hat E}}

\def\oF{{\overline{F}}}
\def\oG{{\overline{G}}}
\def\CC{{\cal C}}
\def\uF{{\underline{F}}}
\def\uG{{\underline{G}}}
\def\oC{{\overline{C}}}
\def\uC{{\underline{C}}}

\begin{document}
\title{On the expected diameter of an $L_2$-bounded martingale}

\setcounter{page}{0}

\baselineskip= 20pt

\author{
Lester E. Dubins \\
Department of Mathematics \\
University of California at Berkeley \\
Berkeley, CA 94720, USA \\
and \\
David Gilat and Isaac Meilijson \\
School of Mathematical Sciences \\
R. and B. Sackler Faculty of Exact Sciences \\
Tel Aviv University, Tel Aviv 69978, Israel \\
{\em
E-mail: \tt{isaco@math.tau.ac.il}} \\
{\em E-mail: \tt{gilat@math.tau.ac.il}}}

\maketitle





\begin{abstract}

\medskip

\begin{center}
{\bf Dedicated to the memory of Gideon Schwarz (1933-2007)}
\end{center}

\medskip

\noindent It is shown that the ratio between the expected diameter
of an
 $L_2$-bounded martingale and the standard deviation of
its last term cannot exceed $\sqrt{3}.$ Moreover, a one-parameter
family of stopping times on standard Brownian Motion is exhibited,
for which the $\sqrt{3}$ upper bound is attained. These stopping
times, one for each cost-rate $c,$ are optimal when the payoff for
stopping at time $t$ is the diameter $D(t)$ obtained up to time $t$
minus the hitherto accumulated cost $c t.$ A quantity related to
diameter, {\em maximal drawdown} (or {\em rise}), is introduced and
its expectation is shown to be bounded by $\sqrt{2}$ times the
standard deviation of the last term of the martingale. These results
complement the Dubins \& Schwarz respective bounds $1$ and
$\sqrt{2}$ for the ratios between the expected maximum and maximal
absolute value of the martingale and the standard deviation of its
last term. Dynamic programming (gambling theory) methods are used
for the proof of optimality.

\medskip



\end{abstract}


\setlength{\textwidth}{6in}
\setlength{\evensidemargin}{0.2in}
\setlength{\oddsidemargin}{0.2in}
\setlength{\topmargin}{0.0in}
\setlength{\textheight}{9in}
\setlength{\headsep}{10pt}
\setlength{\columnsep}{0.375in}
\def\IR{\Bbb R}
\def\tF{{\tilde F}}
\def\tG{{\tilde G}}
\def\tE{{\tilde E}}
\def\hF{{\hat F}}
\def\hG{{\hat G}}
\def\hE{{\hat E}}

\def\oF{{\overline{F}}}
\def\oG{{\overline{G}}}
\def\CC{{\cal C}}
\def\uF{{\underline{F}}}
\def\uG{{\underline{G}}}
\def\oC{{\overline{C}}}
\def\uC{{\underline{C}}}

\section{Introduction}

Lester Dubins \& Gideon Schwarz (\cite{D&S}, 1988) prove that the
ratio between the expectation of the maximum $M$ of a mean-zero
$L_2$-bounded martingale (thus, uniformly integrable, with a well
defined terminal element) and the standard deviation ($L_2$-norm)
of its last term is bounded above by $1.$ They go on to show that
this bound is attained by the martingale $\{B(t) : t \le \tau \},$
where the process $B = \{B(t) : t \ge 0 \}$ is standard Brownian
Motion and $\tau=\tau_d,$ given by
\begin{equation} \label{Tausubd}
\tau_d=\inf\{t \ge 0 : M(t)-B(t)\ge d\} \ , \end{equation} is the
first time $B$ displays a {\em drawdown} of size $d,$ i.e., $B$
drops $d$ units below the highest position it has visited so far;
here $M(t)$ is the maximum of $B$ on $[0,t]$ while $d$ is any
positive constant.

Clearly, the dual stopping time $\tau'$ for minimizing the
expected minimum $m$ relative to the standard deviation of the
last term would be $\tau'=\tau_d'=\inf\{t \ge 0 : B(t)-m(t) \ge d
\}$ , $m(t)$ being the minimum of $B$ on $[0,t].$ $\tau'$ is the
first time $B$ displays a {\em rise} of size $d.$

The main purpose of the present article is to demonstrate that $B$
stopped at time
\begin{equation} \label{calT}
{\cal T} = {\cal T}_d = \inf\{t \ge0 : (M(t)-B(t)\wedge(B(t)-m(t))
\ge d \}
\end{equation}
attains the least upper bound, whose value will be shown to be
$\sqrt{3},$ on the ratio between the expected diameter ($D=M-m$) and
the standard deviation of the last term of any $L_2$-bounded
martingale.

It is useful to point out that the stopping time ${\cal T}$ can be
implemented in two stages: first wait until for the first time a
diameter of size $2d$ is obtained; at this moment $B$ must be
either at its hitherto maximum (up) or minimum (down); if it is
up, continue until from that time on a drawdown of size $d$ is
displayed; similarly, if it is down, continue until a rise of size
$d$ is displayed. It is easy to check that this $2$-stage
procedure terminates exactly at time ${\cal T}.$

Dubins \& Schwarz \cite{D&S} consider also the analogous
inequality for the \linebreak expected supremum $S$ of a
nonnegative $L_2$-bounded submartingale (which by Gilat
(\cite{Gilat1}, 1977) is the same as the absolute value of a
martingale), proving that the least upper bound on the ratio
between $E[S]$ and the
$L_2$-norm (square-root of the
second moment) of the last term is $\sqrt{2}.$
Moreover, as they show, this bound is attained by the absolute
value $|B|$ of $B,$ stopped at time
\begin{equation} \label{Tplus}
T = T_d = \inf\{t \ge 0 : S(t)-|B(t)| \ge d \}
\end{equation}
where $S(t)$ is the supremum of $|B|$ on $[0,t]$ and as before,
$d$ is any positive constant.


Rephrasing $S$ and $D$ in terms of $M$ and $m$ by setting $S=M
\vee |m|$ and $D = M-m = M+|m|,$ it is seen that the stopping time
$\tau$ is optimal for $M,$ its dual $\tau'$ is optimal for $|m|,$
$T$ is optimal for the maximum $M \vee |m|$ of $M$ and $|m|$, and
${\cal T}$ for their sum $M +|m|.$ {\em Optimal} here means
maximizing the pertinent ratio of expectation to standard
deviation. The respective least upper bounds for these ratios are
$1 , 1, \sqrt{2}$ and $\sqrt{3},$ the last being the main
contribution of the present paper.

Related to the diameter $D$ of a process $X$ (with $X(0)=0$) are its
one-sided versions, the {\em maximal drawdown} $D^+$ and the {\em
maximal rise} $D^-$ (with
$D=D^+ \vee D^-$) defined in
terms of $M_X(t)=\sup_{s \le t} X(s)$ and
$m_X(t)=\inf_{s \le t} X(s),$ as follows
\begin{eqnarray}
D^+ & = & \sup_{t \ge 0} \{X(t) - \inf_{s>t} X(s)\} = \sup_{t \ge
0} \{M_X(t) - X(t)\}
\label{Dplusminus} \\
D^- & = & \sup_{t \ge 0} \{\sup_{s>t} X(s) - X(t)\} = \sup_{t \ge
0} \{X(t) -m_X(t)\} \nonumber
\end{eqnarray}

It will be shown (Theorem \ref{optimalstopDplus} \& Corollary
\ref{expecteddrop}) that the supremum, over all $L_2$-bounded
martingales $X,$ of the ratio between $E[D^+]$ (similarly for
$E[D^-]$) and the standard deviation of the last term of $X$ is
$\sqrt{2}.$ Furthermore, this bound is attained by the martingale
$\{B(t) : t \le {\cal T}^+\},$ where
\begin{eqnarray} \label{calTdef}
{\cal T}^+ = {{\cal T}_d}^+ & = & \inf \{t \ge 0 : \sup_{s \le t}
A(s) - A(t) \ge d \} \nonumber \\
& = &  \inf \{t>\tau_d : B(t) - \inf_{\tau_d < s \le t} B(s) \ge d\}
\end{eqnarray}
is the earliest time the {\em drop process} $A(t)=M(t)-B(t)$
attains a drawdown of size $d.$ Equivalently, ${\cal T}^+$ is the
earliest time $B$ attains a rise of size $d$ after having had a
drop of size $d.$

Recalling that the variance of $B(t)$ is $t$ (in fact, $\{B^2(t) - t
: t \ge 0\}$ is a mean-zero martingale), it can be seen (as also
observed in Dubins \& Schwarz \cite{D&S} for maximizing the expected
maximum) that the problem of maximizing the desired ratios is
closely related to that of finding an optimal stopping time on $B$
for the payoff function $ R(t) - ct,$ $c>0$ being the cost per unit
time of sampling. For brevity, refer to this as the $c$-{\em
problem}. Here the {\em reward function} $R(t)$ can be any of the
quantities $M(t),$ $m(t),$ $S(t),$ $D(t)=M(t)-m(t)$ and its two
one-sided versions. In fact, it is not hard to obtain the solution
to the ratio-maximization problem from the solution to the
corresponding $c$-problem and vice versa. We choose to focus on the
latter because it can be conveniently formulated as a continuous
time dynamic programming (or gambling) problem, for which a toolkit
is readily available.

\medskip

\noindent {\bf Remark}. Dubins \& Schwarz \cite{D&S} solve the
ratio-maximization problem for $M$ directly and then infer the
solution to the corresponding $c$-problem. For $S,$ however, they
go only in the opposite direction. A direct solution to the
ratio-maximization for $S$ can be found in Gilat (\cite{Gilat2},
1988). We do not know how to solve the ratio-maximization problem
for $D$ other than by first solving the corresponding $c$-problem.

\medskip

Recalling the definition (\ref{calT}) of the stopping time ${\cal
T},$ our main result now follows:

\begin{theorem} \label{optimalstopD}

$\mbox{    }$

\noindent {\em (}{i}{\em )} For each $c>0,$ ${\cal T} = {\cal
T}_{1 \over {2 c}}$ is optimal for the $c$-problem with reward
function 
$R(t)=D(t).$ It is the unique optimal stopping time within the
family $\{{\cal T}_d : d>0\}.$

\noindent {\em (}{ii}{\em )} $E[{\cal T}] = {3 \over {4c^2}}$

\noindent {\em (}{iii}{\em )} The optimal expected payoff is
$E[D({\cal T}) - c {\cal T}] = {3 \over {4c}}$

\end{theorem}

\begin{corollary} \label{expecteddiam}
The expected diameter of any  $L_2$-bounded martingale cannot
exceed $\sqrt{3}$ times the standard deviation of its last term.
The upper bound  $\sqrt{3}$ is attained by the segment of Brownian
Motion between zero and any of the stopping times ${\cal T}_d.$
\end{corollary}

To state the next Theorem, let $D^+(t)$ be the $D^+$ variable (see
(\ref{Dplusminus})) defined on the martingale $X(\cdot)=B(\cdot
\wedge t),$ and recall the definition (\ref{calTdef}) of the
stopping time ${\cal T}^+.$

\begin{theorem} \label{optimalstopDplus}

$\mbox{    }$

\noindent {\em (}{i}{\em )} For each $c>0,$ ${\cal T^+} = {\cal
T}_{1 \over {2
 c}}^+$ is optimal for the $c$-problem with reward function
$R(t)=D^+(t).$ It is the unique optimal stopping time within the
family $\{{\cal T}_d{^+} : d>0\}.$

\noindent {\em (}{ii}{\em )} $E[{\cal T^+}] = {2 \over {c^2}}$

\noindent {\em (}{iii}{\em )} The optimal expected payoff is
$E[D^+({\cal T^+}) - c {\cal T^+}] = {1 \over {2 c}}$

\end{theorem}

\begin{corollary} \label{expecteddrop}
The expected maximal drawdown of any $L_2$-bounded
martingale cannot exceed $\sqrt{2}$ times the standard
deviation of its last term. The upper bound $\sqrt{2}$ is attained
by the segment of Brownian Motion between zero and any of the
stopping times ${\cal T}_d{^+}.$
\end{corollary}

\noindent {\bf Comment on the special relevance of Brownian
Motion}. It should not be surprising that in a variety of
martingale inequalities (those considered here included), the
extremal martingales, namely those for which equality is attained,
are segments of Brownian Motion determined by a suitable stopping
time. Moreover, in order to establish an inequality for the class
of $L_2$-bounded martingales, it typically suffices to consider
the subclass of these processes of the form $\{B(t):t \leq T\},$
where $T$ is a stopping time with $E[T]< \infty .$ This is so
simply because Brownian Motion is a universal martingale in the
following very specific sense. Recall the Skorokhod (\cite{Skoro},
1965) embedding of a random variable $Z$ with $E[Z]=0$ and
$E[Z^2]< \infty$ in Brownian Motion by a stopping time $T,$ such
that $B(T) \sim Z$ and $E[T]=E[Z^2]$; following I. Monroe
(\cite{Monroe}, 1972) call such a stopping time {\em {minimal}}
(for $Z$). Monroe (\cite{Monroe}, Theorem 11) extends Skorokhod's
result as follows: given a right-continuous, mean-zero,
$L_2$-bounded martingale $X,$ there exists an increasing family
$\{T_t:t \geq 0\}$ of minimal stopping times such that the
embedded process $\{B(T_t)\}$ has the same distribution as $X.$ By
$L_2$-boundedness it follows that the limiting stopping time
$T=\lim_{t \rightarrow \infty} T_t$ is minimal and that $B(T)$ has
the same distribution as the last term of $X.$ Note also that a
process in discrete-time can always be extended to continuous time
and made right-continuous by setting it constant between
consecutive integer time-points. Clearly, the maximum or the
diameter of the entire Brownian path up to time $T,$ dominates the
respective quantities in any embedded process. Consequently, it is
enough to establish our inequalities for Brownian Motion stopped
at minimal stopping times.

\section{Excessivity and supermartingales - Proofs}

\subsection{Proof of Theorem \ref{optimalstopDplus} and its corollary}

Recalling the definition of $D^+(t)$ preceding the statement of
Theorem \ref{optimalstopDplus} and the definition of the drop
process $A$ following (\ref{calTdef}), observe that
$D^+(t)=\sup_{s \le t} A(s).$ Since $A$ is distributed like the
absolute value of a Brownian Motion (see Karatzas \& Shreve
\cite{KarSh} p. 97, who attribute this result to Paul L\'{e}vy),
the $c$-problems for maximal drawdown and maximal absolute value
are equivalent. Consequently, Theorem \ref{optimalstopDplus} and
its corollary follow from the $\sqrt{2}$- inequality of Dubins \&
Schwarz \cite{D&S} (quoted in the Introduction) regarding the
absolute value of a martingale.

\subsection{Proof of Theorem \ref{optimalstopD}}

Define a real-valued function $q = q_{c,d}$ on the domain
$\{(\delta,\gamma,t) : 0 \le \gamma \le {\delta \over 2} < \infty,
t \ge 0\}$ in ${\cal R}^3$ by
\begin{equation} \label{Qfunc}
q(\delta,\gamma,t)= \delta - c t + \left\{
\begin{array}{ll}
0 & \gamma \ge d \\
3 d - \delta - c \{\gamma(\delta-\gamma) + 3 d^2 - {{\delta^2}
\over 2}\}
& \delta < 2 d \\
(d - \gamma)[1 - c(d+\gamma)] & \delta \ge 2 d , \gamma < d
\end{array} \right.
\end{equation}

Note that $q$ is a continuous function.

Let $D(t)=M(t)-m(t)$ be the diameter attained by $B$ by time $t$
and let $G(t)=(M(t)-B(t))\wedge(B(t)-m(t))$ be the {\em gap}, or
minimal distance of the current position from the extremal points
visited so far. Consider the process
$Q(t)=Q_{c,d}(t)=q_{c,d}(D(t),G(t),t)$ and set
$\Pi(t)=D(t)-ct$, the payoff function.

With the help of Lemmas \ref{Qlemma} and \ref{Qdeltalemma} and
Corollary \ref{deltadelta} below, $Q$ can be identified as the
conditional expected payoff for the $c$-problem (with reward
$R(t)=D(t)$) given a partial history $\{B(s) : s \le t\}$ with
current diameter $D(t)$ and gap $G(t),$ when the following stopping
time $\tau_{c,d,t}$ is used: if $G(t) \ge d,$ $\tau_{c,d,t}=t$ ;
otherwise, $\tau_{c,d,t}$ is the first time after $t$ at which the
gap $G$ is at least $d.$ In other words, $\tau_{c,d,t}$ extends
${\cal T}_d$ (see (\ref{calT})) to general initial conditions. That
${\cal T}_d,$ with ${d={1\over{2c}}},$ is optimal for the
$c$-problem will follow from properties of the $Q$ process to be
established in Proposition \ref{mainpropo}: $Q$ majorizes the payoff
$\Pi,$ $Q(0)$ is the expected payoff when using ${\cal T}_d$ and $Q$
is a supermartingale. Thus, for every integrable stopping time
$\tau,$ $E[\Pi(\tau)] \le E[Q(\tau)] \le Q(0) = E[\Pi({\cal T}_d)].$
$Q$ being a supermartingale is the same as $q$ being {\em excessive}
in the Gambling-theoretic terminology of Dubins \& Savage
(\cite{DuSav}, 1965, 1976), a notion closely related to {\em
no-arbitrage pricing} in Finance (see, e.g., pp. 92 in Dana \&
Jeanblanc \cite{DanaJean}, 2003).

\medskip

The following two lemmas summarize known results, most of which are
used in the sequel.

\begin{lemma} \label{Qlemma}
Recall {\em (\ref{Tausubd})} and let $\epsilon_{a,b,x},$ with $a \le
x \le b,$ be the first exit time from the interval $(a,b)$ by
Brownian Motion starting at $x.$

\noindent {\em (}{i}{\em )} (Common knowledge, see e.g.
\cite{Freedman}, p. $71$) $E[\epsilon_{a,b,x}] = (x-a)(b-x)$ and
$P(B(\epsilon_{a,b,x}) = a) = {{b-x} \over {b-a}}$. Similarly for
simple random walk when $a,b$ and $x$ are integers.

\noindent {\em (}{ii}{\em )} (Dubins \& Schwarz {\em \cite{D&S}})
$M(\tau_d)=B(\tau_d)+d$ is exponentially distributed with mean $d.$
Hence, $E[M(\tau_d)]=d$ and $E[\tau_d] = \mbox{Var}[B(\tau_d)]=d^2$.
\end{lemma}

\begin{lemma} \label{Qdeltalemma}
Let $\delta_h$ be the first time $B$ attains a diameter of size $h$.

\noindent {\em (}{i}{\em )} (Pitman \cite{Pitman}). $M(\delta_h)$
and $m(\delta_h)$ are uniformly distributed on their
respective ranges $[0,h]$ and $[-h, 0]$.

\noindent {\em (}{ii}{\em )} (Imhof \cite{Imhof}). The distribution
of $B(\delta_h)$ is given by the V-shaped density function $f_h(x) =
{|x| \over h^2} \ , \ |x| \le h$. Consequently,
$E[\delta_h]=E[B^2(\delta_h)]={h^2 \over 2}$. Similarly, for
positive integer $h$ and integer $x \in [-h, h]$, the probability
that simple random walk stopped at $\delta_h$ terminates at $x$ is
${|x| \over {h(h+1)}}$.
\end{lemma}

Imhof (\cite{Imhof}, formula (2.1)) identifies the joint
distribution of ($\delta_h \ , \ B(\delta_h)$) and obtains (formula
(2.2)) the V-shaped marginal density of $B(\delta_h)$. Here is a
direct argument for random walk, from which the statement for $B$
follows by a standard limiting argument: for $x \in \{1,2,\cdots
h\}$ (similarly for $x \in \{-h, -h+1,\cdots,-1\}$), termination
occurs at $x$ if and only if $x-h$ is reached before $x$ and then
$x$ is reached before going below $x-h$. Since the probability of
the second stage is independent of $x$, the probability of
terminating at $x$ is proportional to the probability of the first
stage, which by Lemma \ref{Qlemma} ({\em i}) is ${x \over
{x+(h-x)}}= {x \over h}$. The consequence for $E[\delta_h]$ is
implied by $B^2(t)-t$ being a mean-zero martingale, and the V-shaped
density having variance ${h^2 \over 2}$.

Pitman (\cite{Pitman} p. 322) infers Lemma \ref{Qdeltalemma} ({\em
i}) from ({\em ii}) in the framework of Brownian Motion on a circle,
when first covering the entire circle. Here is a direct argument:
for $0 < x < h$, $M(\delta_h) \le x$ iff $B$ reaches $x-h$ before
$x$. By Lemma \ref{Qlemma}({\em i}), this event has probability $x
\over h$.

\begin{corollary} \label{deltadelta}
{\em (}{i}{\em )}The expected additional time Brownian Motion needs
to increase its diameter from $h_1$ to $h_2 > h_1$ is
$E[\delta_{h_2}]-E[\delta_{h_1}] = {{h_2^2-h_1^2} \over 2}$

\noindent {\em(}{ii}{\em)} ${E[{\cal
T}_d]=E[\delta_{2d}]+E[\tau_d]=3d^2}$

\noindent {\em(}{iii}{\em)} $E[D({\cal T}_d)] = 3 d$
\end{corollary}

\noindent {\bf Proof}. Claim ({\em i}) follows directly from Lemma
\ref{Qdeltalemma} ({\em ii}). By the two-stage description of ${\cal
T}_d$ which follows its definition (\ref{calT}), ${\cal T}_d$ is the
sum of $\delta_{2d}$ and a random time distributed like $\tau_d$.
Claim ({\em ii}) now follows from Lemma \ref{Qlemma} ({\em ii}) and
\ref{Qdeltalemma} ({\em ii}) by taking expectations. The diameter at
time ${\cal T}_d$ consists of the initial $2 d$ plus the increment
obtained during the second stage. By Lemma \ref{Qlemma} ({\em ii})
this increment has mean $d,$ verifying claim ({\em iii}).

\medskip

The next lemma is instrumental in proving (in the following
Proposition \ref{mainpropo}) that the process $Q_{c , {1 \over {2
c}}}$ is a supermartingale.

\begin{lemma} \label{Hermitelemma}
\noindent {\em (}{i}{\em )} (Paul L\'{e}vy \cite{Levy}) The
processes $\{|B(t)| : t \ge 0\}, \{M(t)-B(t) : t \ge 0\},
\{B(t)-m(t) : t \ge 0\}$ are identically distributed.

\noindent {\em (}{ii}{\em )} The processes $\{B^2(t)-t : t \ge 0\},
\{(M(t)-B(t))^2-t : t \ge 0\}, \{(B(t)-m(t))^2-t : t \ge 0\}$,
adapted to the filtration of $B$, are mean-zero martingales.

\noindent {\em (}{iii}{\em )} The processes $\{\max(B(t),0)^2-t : t
\ge 0\},\{\min(B(t),0)^2-t : t \ge 0\}$ are supermartingales.
\end{lemma}

Proof: Assuming that the martingale nature of $\{B(t)^2 - t\}$ is
well known, statement ({\bf ii}) follows from ({\bf i}). To prove
Lemma \ref{Hermitelemma} ({\em iii}), let $0 \le t < s$. In terms of
the stopping time $\rho=\min(s,\inf \{u : u
> t , B(u) \ge 0\})$, there are three possible cases to consider: $\{B(t) \ge
0\}$, $\{B(t) < 0, \rho=s\}$ and $\{B(t) < 0, t<\rho<s\}$. In the
first case, by statement ({\em ii}) of the Lemma,
$E[\max(B(s),0)^2-s | {\cal B}_t]\le E[B^2(s)-s | {\cal B}_t]=
B^2(t)-t = \max(B(t),0)^2-t$ \linebreak a. s. For the other two
cases, condition on ${\cal B}_{\rho}$ (which contains ${\cal B}_t$)
to obtain
\begin{eqnarray} E[\max(B(s),0)^2-s | {\cal B}_t] & = &
E[E[\max(B(s),0)^2-(s-\rho)|{\cal B}_\rho] | {\cal B}_t] - t
\nonumber
\\
& - & E[\rho - t | {\cal B}_t] \nonumber \\
& < & E[E[\max(B(s),0)^2-(s-\rho)|{\cal B}_\rho] | {\cal B}_t]-t
\label{Hermite}
\end{eqnarray}

In the second case, $\max(B(s),0)^2-(s-\rho)=0=\max(B(t),0)^2$ a. s.
In the third case, by statement ({\em ii}) of the Lemma,
$E[\max(B(s),0)^2-(s-\rho)|{\cal B}_\rho] \le E[B^2(s)-(s-\rho) |
{\cal B}_\rho]=B^2(\rho)=0=\max(B(t),0)^2$ a.s. So in each of these
two cases we obtain that the RHS of (\ref{Hermite}) is bounded from
above by $\max(B(t),0)^2-t$.

\begin{proposition} \label{mainpropo}

For $d={1 \over {2 c}},$ the process $Q= Q_{c,d}=Q_{c,{1 \over {2
c}}}$ has the following properties:

\noindent {\em (}{i}{\em )} $\forall t \ge 0 \ , \ Q(t) \ge
\Pi(t)$ a.s.

\noindent {\em (}{ii}{\em )} $Q(0)=E[\Pi({\cal T}_{1 \over {2 c}}
)]=E[Q({\cal T}_{1 \over {2 c}})].$ Moreover, $\Pi({\cal T}_{1
\over {2 c}})=Q({\cal T}_{1 \over {2 c}})$ a.s.

\noindent {\em (}{iii}{\em )} $Q$ is a supermartingale and
$\{\check{Q}(t)=Q(t \wedge {\cal T}_{1 \over {2 c}}) : t \ge 0\}$
is a martingale (w.r.t. the filtration $\{{\cal B}_t\}$ of the
underlying Brownian Motion)
\end{proposition}


\medskip

\noindent {\bf Proof}. Substituting ${1\over 2c}$ for $d$ in
(\ref{Qfunc}), a straightforward calculation yields
\begin{small}
\begin{eqnarray} \label{Qfunc2c}
& & Q(t)-\Pi(t)= \\
& & \left\{ \begin{array}{ll}
0 & G(t) \ge {1\over 2c} \\
c {[{1 \over 4}({1 \over c}-D(t))({3\over c}-D(t))+({D(t)\over 2}
-G(t))^2]} &
D(t) < {1\over c} \\
c [{1\over 2c} - G(t)]^2
& D(t) \ge {1\over c}
, G(t) < {1\over 2c}
\end{array} \right. \nonumber
\end{eqnarray}
\end{small}

\noindent which is nonnegative, thus verifying claim ({\em {i}}).
Since ${G(0)=D(0)= {\Pi (0)}=0},$ it follows from (\ref{Qfunc2c})
that ${Q(0)={3\over 4c}}.$ The first equality in claim ({\em
{ii}}) now follows by applying ({\em ii}) and ({\em {iii}}) of
Corollary 3 (with ${d={1\over{2c}}}$). To prove the third, a
fortiori the second equality in claim ({\em {ii}}), just note that
by definition $G({\cal T}_{1 \over {2c}})={1 \over {2c}}$ and
$D({\cal T}_{1 \over{2c}}) \ge {1 \over c}$ .

To establish claim (iii), the time-axis $[0,\infty)$ will be
partitioned into a sequence of intervals with suitably chosen
stopping times as their end-points. The process $Q(\cdot)$ will then
be represented in an appropriate form, tailored for the application
of Lemma \ref{Hermitelemma}, over each of these subintervals. To
exhibit this partition, fix an arbitrary $f \in (0,d)$ and
inductively define an increasing sequence of stopping times as
follows: $\tau_0=0$ a.s., $\tau_1$ is the first time $B$ achieves
diameter $2 d$ and $\tau_2 ={\cal T}_d$. Note that $B(\tau_2)$ is an
endpoint of the a.s. non-empty {\em central interval}
$(m(\tau_2)+d,M(\tau_2)-d)$. Let $\tau_3$ be the earlier between the
next time $B$ reaches the other endpoint of the central interval or
the gap decreases to $f$. Generically now, for $n \ge 3$, if
$B(\tau_{n-1})$ and $B(\tau_{n})$ are the endpoints of the central
interval $(m(\tau_n)+d,M(\tau_n)-d)$, define $\tau_{n+1}$ similarly
to $\tau_3$. If on the other hand the gap at $B(\tau_n)$ is $f$, let
$\tau_{n+1}$ be the first time after $\tau_n$ at which the gap
reaches $d$ again. We now use (\ref{Qfunc2c}) to represent
$Q(\cdot)$ over each of the partition intervals
$[\tau_{n-1},\tau_n)$ in a form conducive to the application of
Lemma \ref{Hermitelemma}.

Between times $\tau_0$ and $\tau_1$, $Q$ is equal to the martingale
(see Lemma \ref{Hermitelemma} ({\em ii})) ${3 \over {4 c}} + {c
\over 2} \{ [(M(t)-B(t))^2 - t] + [(B(t)-m(t))^2 - t] \}$. Resorting
to the short-hand $t^*=\max(\tau_1,\min(\tau_2,t))$, between times
$\tau_1$ and $\tau_2$, $Q$ is equal to the martingale defined as $Q$
up to $\tau_1$ and, thereafter (see again Lemma \ref{Hermitelemma}
({\em ii})),
${1 \over {4 c}} + B(t^*) - m(\tau_1) + c [(M(t^*) - B(t^*))^2 -
t^*]$ or its mirror image, depending on $B(\tau_1)$ being
$M(\tau_1)$ or $m(\tau_1)$.
A similar representation of $Q(\cdot)$ is readily available for the
time-increments during which the gap increases from $f$ to $d$.
Finally, consider a time $\tau_n$ at which $B$ is at an endpoint of
the then central interval. Letting
$t^*=\max(\tau_n,\min(\tau_{n+1},t))$, from time $\tau_n$ to time
$\tau_{n+1}$, the process $Q$ is equal to the supermartingale
%
defined as $Q$ up to $\tau_n$, and thereafter (see Lemma
\ref{Hermitelemma} ({\em iii})), $M(\tau_n) - m(\tau_n) - c \tau_n +
c [ \max(W(t^*),0)^2 - (t^* - \tau_n)]$, where, depending on
$B(\tau_n)$ being $M(\tau_n)+d$ or $m(\tau_n)-d$, $W(\cdot) =
B(\cdot) - B(\tau_n)$ or its mirror image.

Since $\cup_{n=1}^\infty [\tau_{n-1},\tau_n) = [0,\infty)$ a.s.,
$Q(\cdot)$ is a supermartingale throughout.

\medskip

\noindent {\bf Concluding the proof of Theorem \ref{optimalstopD}}.
As argued prior to the statement of Lemma \ref{Qlemma}, Proposition
\ref{mainpropo} establishes the optimality of ${\cal T}_{1 \over {2
c}}$ for the current $c$-problem.

By Corollary \ref{deltadelta}, $E[\Pi({\cal T}_d)]=E[D({\cal
T}_d)]-c E[{\cal T}_d]=(2 d + d)-c({{(2 d)^2} \over 2}+d^2)=3 d - 3
c d^2,$ which is uniquely maximized at $d={1 \over {2 c}}.$ Thus
${\cal T}_{1 \over {2 c}}$ is the unique optimal stopping time
within the family $\{{\cal T}_d \ : \ d>0\}.$

\noindent Claims ({\em ii}) and ({\em iii}) of Theorem
\ref{optimalstopD} follow straightforwardly from Lemma
\ref{Qdeltalemma}.

\medskip

\noindent {\bf Proof of Corollary \ref{expecteddiam}}. As argued
in the last paragraph of the Introduction, it is enough to
consider martingales of the form $\{B(t) : t \le T\},$ where $T$
is a stopping time with $E[T]=E[(B(T))^2] < \infty.$ Let
$\sigma=\sqrt{E[T]}$ and consider the $c$-problem with
$c={\sqrt{3} \over {2 \sigma}}.$ Then
\begin{eqnarray} \label{EDT}
E[D(T)] & = & (E[D(T)]-c E[T])+c E[T] \nonumber \\
& \le & E[D({\cal T}_{1 \over {2 c}})] - c E[{\cal
T}_{1 \over {2 c}}]+c \sigma^2 \nonumber \\
& = & {3 \over {4 c}}+c \sigma^2 = {3 \over {4 {\sqrt{3} \over {2
\sigma}}}} + {\sqrt{3} \over 2} \sigma = \sqrt{3} \sigma
\end{eqnarray}
That the $\sqrt{3}$-bound is attained by$\{B(t) : t \le {\cal
T}_d\}$ for any $d>0$ follows from Corollary \ref{deltadelta}.

\section{An open problem: the spider process}

Larry Shepp has recently reminded us that the $\sqrt{3}$-inequality
treated here is a special case of the so called {\em spider} problem
raised some time ago by the first author. Informally speaking,
Brownian Motion may be viewed as an absolute value of Brownian
Motion, each of whose excursions is assigned a random sign. The
spider process (sometimes called the {\em Walsh} process) with $n
\ge 3$ rays emanating from the origin is the extension from Brownian
Motion ($n=2$) to an $n$-valued sign. (Thus, $n=4$ corresponds to
randomly switching at visits to zero between an absolute value of BM
on the y-axis and another on the x-axis). The maximal distance from
the origin in the spider process is simply the maximal absolute
value of Brownian Motion, independently of $n.$ On the other hand,
the sum of the distances from the origin along the rays, reduces in
the case $n=2$ to the diameter of Brownian Motion studied here. The
maximization of the expected value of this sum of distances when $n
\ge 3$ seems harder to handle and evidently requires new ideas.

\section*{Acknowledgements}
We thank Jim Pitman for helping us avoid some local martingale
traps, Greg Lawler and Larry Shepp for informative
discussions, Heinrich von Weizsaecker for pointing out a minor
mistake in an earlier draft of the introductory section and two
anonymous referees for a thorough review of the paper and useful
suggestions for improvement.


\end{document}